# ASYMPTOTIC LAWS FOR COMPOSITIONS DERIVED FROM TRANSFORMED SUBORDINATORS[1]

### By Alexander Gnedin, Jim Pitman and Marc Yor

*Utrecht University, University of California and University of Paris VI*

A random composition of $n$ appears when the points of a random closed set $\widetilde{\mathcal{R}} \subset [0,1]$ are used to separate into blocks $n$ points sampled from the uniform distribution. We study the number of parts $K_n$ of this composition and other related functionals under the assumption that $\widetilde{\mathcal{R}} = \phi(S_\bullet)$, where $(S_t, t \geq 0)$ is a subordinator and $\phi : [0, \infty) \to [0,1]$ is a diffeomorphism. We derive the asymptotics of $K_n$ when the Lévy measure of the subordinator is regularly varying at 0 with positive index. Specializing to the case of exponential function $\phi(x) = 1 - e^{-x}$, we establish a connection between the asymptotics of $K_n$ and the exponential functional of the subordinator.

**1. Introduction.** A *composition* of $n$ with positive integer parts may be represented by a configuration of stars separated by bars, for instance, $***|**|****|*$ encodes the composition $(3, 2, 4, 1)$ with *weight* 10, *length* 4 and four *parts* 3, 2, 4, 1. A stochastic analogue of this construction appears when we assume the points of a closed random set $\widetilde{\mathcal{R}} \subset [0,1]$ in the role of bars, and $n$ independent random points sampled from the uniform distribution on $[0,1]$ in the role of stars, see [11, 13, 15, 16, 25]. Given this data, we define an ordered partition of the set $\{u_1, \dots, u_n\}$ by assigning two points $u_i < u_j$ to the same block if and only if $u_i$ and $u_j$ are not separated by $\widetilde{\mathcal{R}}$, meaning that $\widetilde{\mathcal{R}} \cap [u_i, u_j] = \varnothing$. That is to say, $u_j \in \mathcal{R}$ forms a singleton block, and if $u_i$ and $u_j$ fall in the same *gap* (open interval component of $[0,1] \setminus \widetilde{\mathcal{R}}$), then these points are assigned to the same block. A composition $\mathcal{C}_n$ is defined to be the record of block sizes, ordered from left to right. Exchangeability in the infinite sample $u_1, u_2, \dots$ results in a simple consistency condition of $\mathcal{C}_n$'s as $n$ varies, that is, the sequence $(\mathcal{C}_n)$ is a *composition structure* in the sense of [11, 12, 16].

Received March 2004; revised December 2004.

[1]Supported in part by NSF Grant DMS-00-71448.

The model just described offers a general framework for a wide range of "species sampling" problems, as studied in statistics and population genetics. In these applications one postulates some idealized infinite population, randomly partitioned into various species, with a total order on the set of species. A sample from such a population is understood as an exchangeable sequence of random variables $(X_j)$, and a composition $\mathcal{C}_n$ is defined as the record of multiplicities of distinct values represented among $X_1, \ldots, X_n$, in the order of increase of the values. Then $(\mathcal{C}_n)$ is a composition structure and by a de Finetti-type result [11], it can be uniquely associated with some random closed set $\widetilde{\mathcal{R}} \subset [0, 1]$, which appears as a way to uniformize the limiting empirical distribution of $(X_j)$.

Let $K_n$ be the length of $\mathcal{C}_n$ (this variable may be interpreted as the number of distinct species in a sample). The growth properties of moments of $K_n$ are sensitive functions of the random set $\widetilde{\mathcal{R}}$. Logarithmic and power-like asymptotics of the moments are known in the case when $\widetilde{\mathcal{R}}$ is derived by scaling the range of a *subordinator* $(S_t, t \in [0, T])$, that is, increasing process with stationary independent increments restricted to a finite time interval [1, 24, 25, 26]. (See [3] for general background on subordinators.)

In this paper we study asymptotic properties of $K_n$ for the random sets obtained by transforming the unrestricted range of a subordinator. Specifically, we consider $\widetilde{\mathcal{R}} = \phi(S_\bullet)$, where $(S_t, t \geq 0)$ is a drift-free subordinator and $\phi \colon [0, \infty] \to [0, 1]$ is a diffeomorhism. We assume the Lévy measure $\nu$ of the subordinator to be regularly varying in the sense that

$$(1) \qquad\qquad \nu[y, \infty] = \ell(1/y) y^{-\alpha}, \qquad y \downarrow 0,$$

where $0 < \alpha \leq 1$ and the function $\ell$ is slowly varying at $\infty$. We also consider the process $K_n(t)$, the number of parts of the partial composition produced by the transformed subordinator restricted to the time interval $[0, t]$. Other quantities of interest are $K_{n,r}$ and $K_{n,r}(t)$, defined as the multiplicity of part $r$ in $\mathcal{C}_n$ and multiplicity in the partial composition, respectively.

We show that, as $n \to \infty$, the length $K_n$ is asymptotic to a power-like, regularly varying function of $n$ multiplied by a random factor $L$. The factor $L$ is identified explicitly as an integral functional of the subordinator. Similar results also hold for $K_{n,r}, K_n(t)$ and $K_{n,r}(t)$. The appearance of a random factor is due to variability in the gap sizes, as can be compared with a result by Karlin [20] which states that the number of distinct values in a large sample from arbitrary nonrandom discrete distribution is asymptotic to the mean number of such values.

In the special case $\phi(y) = 1 - e^{-y}$, the set $\widetilde{\mathcal{R}}$ is a *multiplicative subordinator* and the composition $\mathcal{C}_n$ inherits a characteristic regenerative property from this set [14, 16]. We show that $L$ specializes in this case as the well-known exponential functional of subordinator. The distribution of $L$ is then uniquely determined by the power moments which are given by a known formula reproved here by elementary tools in the case of subordinators.

In the regenerative case, the distribution of $K_n$ is well known for the composition described by Ewens' sampling formula, in which case $K_n$ is of logarithmic growth [1, 25]. More generally, Gnedin [13] has previously shown that the logarithmic growth of $K_n$ is typical when the Lévy measure is finite. For compositions belonging to the two-parameter family [16, 21, 25], the proper format for $K_n$ is $n^\alpha$ for parameters $(\alpha, \theta)$ with $0 < \alpha < 1$. Another interesting case is that of slow variation, when the relation (1) holds with $\alpha = 0$ and some $\ell(1/y)$ exploding at 0. This includes the gamma subordinators whose Lévy measure has a logarithmic singularity. This case is very different from the case of regular variation with positive index $\alpha$ and is being treated separately [2, 17].

We shall be assuming throughout that (1) holds, which entails that the Lévy measure is infinite. When the Lévy parameters $(\nu, \mathtt{d})$ are multiplied by a positive factor $c$, the variables $K_{n,r}$ remain unchanged, but $K_{n,r}(t)$ should be replaced by $K_{n,r}(t/c)$. Basically, we assume that the Lévy measure satisfies $\nu\{\infty\} = 0$ and that the drift coefficient is 0, unless explicitly stated.

It should be mentioned that there are many other constructions of random compositions, but typically these compositions are not consistent as $n$ varies. One obvious possibility, in terms of the "stars and bars" representation, is to exploit the Bernoulli scheme, that is, to allocate a bar at each possible position with fixed probability $p$. (The particular choice $p = 1/2$ corresponds to the uniform distribution on the set of all compositions of weight $n$, see [18].) The expected length of such composition grows linearly with $n$, while, for composition structures, we have $\mathbb{E}K_n = o(n)$, provided the Lebesgue measure of $\widetilde{\mathcal{R}}$ is 0 (which means that the positive frequencies of distinct species sum to 1). See [15] for a complete characterization of the composition structures obtained by truncating a single infinite sequence of stars separated by bars at positions visited by some increasing random process on integers.

The rest of the paper is organized as follows. In the next three sections we modify Karlin's results on occupancy problems, we provide some analysis of the gap counts necessary to apply these results to the composition derived by a general transform of subordinator and we formulate the strong laws for $K_n$ and the like. We specialize then to multiplicative subordinators in Section 5. In Section 6 we continue to consider the regenerative case, but replace fixed-$n$ sample by a Poisson point process, we then analyze recursions for the moments of the length of poissonised composition and show the convergence of the scaled moments of $K_n$.

**2. General strong laws.** Karlin [20] studied the number of different types represented in a sample from a fixed discrete distribution with infinitely many positive masses. His results open a clear path to the strong laws for $K_n$ and $K_{n,r}$. Let $\widetilde{\mathcal{R}}$ be an arbitrary closed subset of $[0,1]$ with zero Lebesgue measure. Let $\mathcal{C}_n$ be the composition derived from $\widetilde{\mathcal{R}}$ by separating uniform points. Conditionally given $\widetilde{\mathcal{R}}$, the number of parts of $\mathcal{C}_n$ is the same as the number of different types represented in a sample from the discrete distribution with masses equal to the gap-sizes. Therefore, by [20], Theorem 8, as $n \to \infty$,

$$(2) \qquad K_n \sim \mathbb{E}(K_n | \widetilde{\mathcal{R}}), \qquad K_{n,r} \sim \mathbb{E}(K_{n,r} | \widetilde{\mathcal{R}}), \qquad r \geq 1,$$

where $\sim$ means that the ratio converges to 1 almost surely. For $x > 0$, let $\widetilde{N}_x$ be the number of gaps of $\widetilde{\mathcal{R}}$ of size at least $x$. The following is a variation of [20], Theorem 1, equation (23) and page 396. See also [25], Lemma 34.

THEOREM 2.1. Let $\ell$ be a positive slowly varying function and $L$ a nonnegative random variable. The convergence

$$\frac{\widetilde{N}_x}{\ell(1/x)x^{-\alpha}} \to L \qquad a.s., x \downarrow 0,$$

with $0 < \alpha < 1$ implies, for $n \to \infty$,

$$\frac{K_n}{n^\alpha \ell(n)} \to \Gamma(1-\alpha)L, \qquad \frac{K_{n,r}}{n^\alpha \ell(n)} \to \frac{\alpha \Gamma(r-\alpha)}{r!}L$$

almost surely, and the same convergence with $\alpha = 1$ implies

$$\frac{K_n}{n\ell^*(n)} \to L, \qquad \frac{K_{n,1}}{n\ell^*(n)} \to L, \qquad \frac{K_{n,r}}{n\ell(n)} \to \frac{1}{r(r-1)}L \qquad \text{for } r > 1$$

almost surely. Here $\ell^*$ is another function of slow variation at $\infty$, defined by the converging integral

$$(3) \qquad \ell^*(t) = \int_0^\infty \frac{e^{-1/y}}{y}\ell(ty)\,dy.$$

PROOF. Let us start with $K_n$. By (2), it is sufficient to determine the asymptotics of conditional expectation. To this end, we introduce a random measure $\gamma$ on $\,]0,1]$ by defining its tail

$$\tilde{\gamma}(x) := \gamma[x,1] = \widetilde{N}_x, \qquad x \in \,]0,1],$$

to be the number of gaps of $\widetilde{\mathcal{R}}$ of size at least $x$. The measure $\gamma$ is atomic and assigns to each $x \in \,]0,1]$ an integer weight equal to the number of gaps of $\widetilde{\mathcal{R}}$ of length $x$. For a particular gap of length $x$, the probability that at least one of $n$ uniform sample points hits this gap is $1 - (1-x)^n$, so

$$(4) \qquad \mathbb{E}(K_n | \widetilde{\mathcal{R}}) = \int_0^1 (1 - (1-x)^n)\gamma(dx) = n \int_0^1 (1-x)^{n-1}\tilde{\gamma}(x)\,dx,$$

where the second equality is obtained by integration by parts. Observe that the formula

$$(5) \qquad \int_0^1 \bar{\bar{\gamma}}(x)\,dx = \int_0^1 x\gamma(dx) = 1$$

simply says that the total length of gaps equals 1, thus, the measure $\bar{\bar{\gamma}}(x)\,dx$, $x \in [0,1]$, is a probability measure with nonincreasing density, which takes only nonnegative integer values. In the last integral in (4) we recognize a Mellin transform and standard Abel–Tauberian arguments (see Appendix) imply that, for $0 < \alpha < 1$,

$$\bar{\bar{\gamma}}(x) \sim x^{-\alpha}\ell(1/x)L \qquad \text{for } x \downarrow 0 \qquad \text{iff}$$

$$n\int_0^1 (1-x)^{n-1}\bar{\bar{\gamma}}(x)\,dx \sim \Gamma(1-\alpha)n^\alpha \ell(n)L \qquad \text{for } n \to \infty$$

and the result follows in this case. In the case $\alpha = 1$, the Mellin integral is asymptotic to the Laplace integral

$$\int_0^\infty e^{-nx}x^{-1}\ell(1/x)\,dx,$$

which converges due to (5), and becomes (3) upon substituting $nx = 1/y$. The slow variation claim for $\ell^*$ is Lemma 4 in [20].

For $K_{n,r}$, we have a similar integral representation

$$(6) \qquad \mathbb{E}(K_{n,r}|\widetilde{\mathcal{R}}) = \binom{n}{r}\int_0^1 x^r(1-x)^{n-r}\gamma(dx),$$

which is obtained by a formal binomial expansion of $1 - (1-x)^n$. The formula follows by observing that a gap of length $x$ is hit by exactly $r$ sample points with probability $\binom{n}{r}x^r(1-x)^{n-r}$. A Tauberian argument applied to the measure $x^r\gamma(dx)$ yields

$$\mathbb{E}(K_{n,r}|\widetilde{\mathcal{R}}) \sim (-1)^{r-1}\binom{\alpha}{r}\mathbb{E}(K_n|\widetilde{\mathcal{R}}),$$

which ends the proof. $\square$

REMARKS. For the two slowly varying functions in the theorem, we have $\ell(t)^*/\ell(t) \to \infty$ as $t \to \infty$ (formula (13) in [20] is misprinted), that is, for $\ell(t) = (\log t)^{-u}$, $u > 1$, we have $\ell^* \sim (u-1)^{-1}(\log t)^{1-u}$. See [6], Chapter 3, for results involving two slowly varying functions like $\ell$ and $\ell^*$. The relation between asymptotics of (4) and (6) is an instance of "smooth variation" properties [6], Section 1.8 of the Bernstein function defined by (4).

Under the assumption of Theorem 2.1, the conditional distribution of $K_n$ given $\widetilde{\mathcal{R}}$ approaches a normal distribution as $n \to \infty$, by [20], Theorem 4 (also see [8], Theorem 2). Karlin's results also imply a multivariate normal limit for the conditional distribution of the sequence $(K_{n,r}, r \geq 1)$.

**3. Counting the gaps.** Let $(S_t, t \geq 0)$ be a subordinator with the drift coefficient $\mathtt{d} = 0$ and a Lévy measure $\nu$ satisfying $\nu\{\infty\} = 0$. The jumps of $(S_t)$ correspond to the gaps of $\mathcal{R} \subset [0, \infty]$, which is a topological Cantor set provided the Lévy measure is infinite.

Let $\phi \colon [0, \infty[ \to [0, 1[$ be a diffeomorphism, that is, a continuously differentiable function satisfying $\phi(0) = 0$, $\phi(\infty-) = 1$ and $\phi'(t) > 0$. For $\widetilde{\mathcal{R}} := \phi(\mathcal{R}) \subset [0, 1]$, the gaps comprising $\widetilde{\mathcal{R}}^c = [0, 1] \setminus \widetilde{\mathcal{R}}$ correspond to the jumps of the subordinator transformed by $\phi$. Let $\widetilde{N}_x(t)$ be the number of jumps of size at least $x$ for the transformed subordinator restricted to $[0, t]$, and let $\widetilde{N}_x = \widetilde{N}_x(\infty)$ be the number of such gaps of $\widetilde{\mathcal{R}}$ without restriction. We are interested in the asymptotics of these gap counts for small $x$. A similar analysis has appeared in [23] in the case of stable subordinators.

The analogous question for the original subordinator $(S_t)$ is easy. Let $N_y(t)$ be the number of gaps of $\mathcal{R}$ of size at least $y$, generated by the subordinator restricted to $[0, t]$. The counting process $(N_y(t), t \geq 0)$ is a Poisson process with rate $\vec{\nu}(y)$, thus, for small $y$, the behavior of this process is ruled by the strong law of large numbers:

$$(7) \qquad N_y(t) \sim \vec{\nu}(y) t, \qquad y \downarrow 0,$$

almost surely for all $t$. We shall see that translating this behavior into similar results for $\widetilde{N}_x(t)$ and $\widetilde{N}_x = \widetilde{N}_x(\infty)$ amounts to a change of variable formula which was stated in [23], albeit under different assumptions on $\phi$.

Speaking more broadly, we may wonder about the conditions on $\phi$ and $\nu$ which imply an asymptotic relation analogous to (7) of the type

$$(8) \qquad \widetilde{N}_x(t) \sim \psi(x) L(t), \qquad x \downarrow 0,$$

where $\psi$ is a scaling function and $(L(t), t \in [0, \infty])$ is a positive random process. A principal new effect appearing in (8), as compared to (7), is that a nonlinear transformation of subordinator leads to a genuinely random scaling limit. The next question to ask is whether such a relation holds with some $L$ for $t = \infty$ and whether $L(\infty-) = L$, and we shall find conditions when this is true.

For a Lévy measure as in (1), we shall use the scaling function

$$\psi(x) = x^{-\alpha} \ell(1/x).$$

3.1. *Finite $t$ formula.* Let $N_y(t_1, t_2)$ be the number of jumps of $(S_t, t \in [t_1, t_2])$ of size at least $y$, with the convention that this is zero for $t_1 > t_2$.

THEOREM 3.1. *If the Lévy measure satisfies* (1), *then for each diffeomorphism* $\phi \colon [0, \infty[ \to [0, 1[$ *and* $0 < t < \infty$, *the convergence, as* $x \to 0$,

$$(9) \qquad \widetilde{N}_x(t)/\psi(x) \to \int_0^t (\phi'(S_u))^\alpha \, du$$

*holds in the mean for each t, as well as almost surely, uniformly in t bounded away from ∞.*

Proof. Consider a partition of $[0, t]$ by points $0 = t_0 < t_1 < \cdots < t_k = t$; with probability 1, each $t_j$ is a continuity point of the subordinator. As is easily seen,

$$N_{x/\phi'(S_{\eta_j})}(t_j, t_{j+1}) \leq \widetilde{N}_x(t_j, t_{j+1}) \leq N_{x/\phi'(S_{\xi_j})}(t_j, t_{j+1}),$$

where $\xi_j$ and $\eta_j$ are the points where $\phi'$ attains the maximum and the minimum on $[S_{t_j}, S_{t_{j+1}}]$, respectively. Taken together with (7), this implies

$$(10) \quad (t_{j+1} - t_j)\bar{\nu}\left(\frac{x}{\phi'(S_{\eta_j})}\right) \ll \widetilde{N}_x(t_j, t_{j+1}) \ll (t_{j+1} - t_j)\bar{\nu}\left(\frac{x}{\phi'(S_{\xi_j})}\right),$$

where the notation $X \ll Y$ for positive random quantities depending on $x$ means that $\operatorname{ess\,sup} X/Y \leq 1$ for $x \downarrow 0$. From this and the assumption on $\nu$,

$$\psi(x) \sum_{j=0}^{k-1} (\phi'(S_{\eta_j}))^\alpha (t_{j+1} - t_j) \ll \widetilde{N}_x(t) \ll \psi(x) \sum_{j=0}^{k-1} (\phi'(S_{\xi_j}))^\alpha (t_{j+1} - t_j).$$

We see that $\widetilde{N}_x(t)$, for $x \downarrow 0$, can be squeezed between an upper and a lower Riemann sum; thus, sending the diameter of partition to zero and using the continuity, we obtain the almost sure convergence (9). Using the obvious bound $\widetilde{N}_x(t) \leq N_{x/\max \phi'}(t)$ where the maximum is taken over $[0, t]$, the convergence in mean follows by dominated convergence. □

There is a minor generalization of the formula for integrals with random upper bound. By *stopping time* $\tau$, we understand a random variable taking values in $[0, \infty]$ and such that $((\tau \wedge t), t \geq 0)$ is adapted to the natural filtration of the subordinator.

Corollary 3.2. *If the Lévy measure satisfies (1), then for each diffeomorphism $\phi : [0, \infty[ \to [0, 1[$ with $\sup \phi' < \infty$, and each stopping time $\tau$ with $\mathbb{E}\tau < \infty$, the convergence, as $x \to 0$,*

$$\widetilde{N}_x(\tau)/\psi(x) \to \int_0^\tau (\phi'(S_u))^\alpha \, du$$

*holds almost surely and in the mean.*

Proof. The almost sure convergence follows from Theorem 3.1. The convergence in the mean is a consequence of

$$\lim_{t \to \infty} \limsup_{x \downarrow 0} \mathbb{E}\widetilde{N}_x(t, \tau)/\psi(x) = 0,$$

which, in turn, follows from

$$\mathbb{E}\widetilde{N}_x(t,\tau) < \mathbb{E}N_{x/\sup\phi'}(t,\tau) = \mathbb{E}N_{x/\sup\phi'}(\tau-t)_+ = \vec{\nu}\left(\frac{x}{\sup\phi'}\right)\mathbb{E}(\tau-t)_+$$

by application of Wald's identity,

$$\mathbb{E}N_y(\tau) = \vec{\nu}(y)\mathbb{E}\tau,$$

and the fact that $\mathbb{E}(\tau-t)_+ \to 0$ for $t \to \infty$.  $\square$

The corollary can be applied to a subordinator killed at an independent time $\tau$. For example, when $\nu\{\infty\} > 0$, the subordinator jumps to infinity at an independent exponential time.

3.2. *Full range formula.* We turn next to finding some conditions for the convergence

$$\widetilde{N}_x/\psi(x) \to \int_0^\infty (\phi'(S_u))^\alpha \, du, \tag{11}$$

in which case (8) holds for $t = \infty$ with $L = \lim_{t\uparrow\infty} L(t)$ given by the integral in (11). One condition which seems very natural is the integrability

$$\mathbb{E}\int_0^\infty (\phi'(S_t))^\alpha \, dt < \infty. \tag{12}$$

Granted this integrability condition only, we failed to prove or disprove whether the convergence holds in full generality, nor is there an obvious sufficient condition which would cover the cases of interest including slowly varying functions $|\log y|$ or $1/|\log y|$ and diffeomorphisms with exponentially decaying or power-like tails.

To secure the convergence, we shall make some additional assumptions about $\phi$ and $\nu$. The analysis is largely simplified by further assuming that the derivative $\phi'$ is a decreasing function; in this case, in (10), we can set $\xi_j = t_j$ and $\eta_j = t_{j+1}$, and then, for any $t_1 < t_2$, conditioning on $S_{t_1}$ yields the inequality

$$\mathbb{E}\widetilde{N}_x(t_1,t_2) < (t_2-t_1)\mathbb{E}\nu\left(\frac{x}{\phi'(S_{t_1})}\right), \tag{13}$$

which is valid for each $x > 0$.

THEOREM 3.3. *Assume that the Lévy measure is as in* (1), *the diffeomorphism* $\phi$ *has decreasing derivative, the integrability condition* (12) *holds, and one of the following single or composite conditions is satisfied:*

(i) *for some constants* $a > 0, C > 0$, *the inequality* $\ell(u) < C\ell(v)$ *holds for* $u < av$, *provided* $u$ *and* $v$ *are sufficiently large;*

(ii) *there exist functions $q$ and $r$ such that:*

(iia) $q(y) = o(\psi(y))$ *as $y \downarrow 0$;*

(iib) *for some constants $a > 0, C > 0$, the inequalities $r(1/v)v < u < av$ imply $\ell(u) < C\ell(v)$ for all sufficiently large $u, v$;*

(iic) $\bar{\nu}(x/r(x))\phi^{\leftarrow}(1-q(y)) = o(\psi(y))$ *as $y \downarrow 0$;*

(iii) $\ell$ *is bounded away from 0 and $\infty$ on every compact subset of $[0, \infty[$ and a stronger integrability condition holds, with $\alpha$ in* (12) *replaced by $\alpha - \delta$, for some $\delta > 0$;*

(iv) *the same integrability condition holds as in* (iii) *and there exists a function $q$ which satisfies* (iia), *as well as $\phi^{\leftarrow}(1-q(y)) = o(\psi(y))$ for $y \downarrow 0$.*

*Then the convergence* (11) *holds almost surely and in the mean.*

PROOF. In view of $\liminf \tilde{N}_x/\psi(x) \geq L(\infty-)$, it is sufficient to establish the convergence of expectations

$$\lim_{x \to 0} \mathbb{E}\tilde{N}_x/\psi(x) = \mathbb{E}\int_0^\infty (\phi'(S_t))^\alpha \, dt.$$

Thanks to both convergence results,

$$\lim_{t \to \infty} \mathbb{E}\int_0^t (\phi'(S_u))^\alpha \, du = \mathbb{E}\int_0^\infty (\phi'(S_u))^\alpha \, du < \infty,$$

$$\lim_{x \to 0} \mathbb{E}\tilde{N}_x(t)/\psi(x) = \mathbb{E}\int_0^t (\phi'(S_u))^\alpha \, du < \infty,$$

and in view of $\tilde{N}_x = \tilde{N}_x(0,t) + \tilde{N}_x(t,\infty)$, we only need to show that

$$\lim_{t \to \infty} \limsup_{x \to 0} \mathbb{E}\tilde{N}_x(t,\infty)/\psi(x) = 0.$$

Using monotonicity,

$$\mathbb{E}\tilde{N}_x(t,\infty) = \sum_{j=0}^\infty \mathbb{E}\tilde{N}_x(t+j, t+j+1) < \sum_{j=0}^\infty \mathbb{E}\bar{\nu}\left(\frac{x}{\phi'(S_{t+j})}\right)$$

$$= \mathbb{E}\sum_{j=0}^\infty (\phi'(S_{t+j}))^\alpha \ell(\phi'(S_{t+j})/x)x^{-\alpha},$$

whence

$$(14) \qquad \mathbb{E}\frac{\tilde{N}_x(t,\infty)}{x^{-\alpha}\ell(1/x)} < \mathbb{E}\sum_{j=0}^\infty (\phi'(S_{t+j}))^\alpha \left(\frac{\ell(\phi'(S_{t+j})/x)}{\ell(1/x)}\right).$$

The rest of the proof amounts to estimating the right-hand side of this formula.

Let $\tau_x$ be the first passage time over the level $(\phi')^{\leftarrow}(x)$, when $\phi'(S_t)$ drops below $x$. As a consequence of Corollary 3.2, the sum of terms in (14) with $t + j \leq \tau_a$ is negligible, for each fixed $a > 0$. And assuming (i), we estimate the contribution of terms in (14) with $t + j > \tau_a$ by

$$C\mathbb{E}\int_{t-1}^{\infty}(\phi'(S_u))^\alpha\,du,$$

which vanishes for $t \to \infty$ due to (12).

Assume (ii). Let $\sigma_z$ be the first passage time over the level $\phi^{\leftarrow}(z)$. By (iia), the contribution of terms with $t + j > \sigma_{1-q(x)}$ is negligible because at most $q(x)/x$ gaps longer than $x$ can fit into an interval of size $q(x)$. The contribution of terms with $t + j \leq \tau_{r(x)}$ also vanishes, as $x \downarrow 0$ and then $t \to \infty$, for the same reason as under the condition (i). Finally,

$$\mathbb{E}\widetilde{N}_x(\tau_x,\sigma_x) < \mathbb{E}\left(\vec{\nu}\left(\frac{x}{\phi'(S_{\tau_x})}\right)(\sigma_x - \tau_x)_+\right)$$

$$< \vec{\nu}\left(\frac{x}{r(x)}\right)\mathbb{E}\sigma_{1-q(x)} = o(x^{-\alpha}\ell(1/x))$$

by virtue of (iic), because the expected time for the subordinator to pass a high level $y = \phi^{\leftarrow}(1 - q(x))$ can be estimated from the above by $by$, with some positive constant $b$.

Suppose (iii) holds. By Potter's Theorem A.3(ii), there exists $C > 1$ such that $\ell(u)/\ell(v) < C(u/v)^{-\delta}$ for all $u < v$. We apply this with $u = \phi'(S_t)/x$ and $v = 1/x$ and then use the integrability with the exponent $\alpha - \delta$.

Under assumption (iv), we make use of another part of the same Theorem A.3(i), which guarantees the same inequality for sufficiently large parameters, say, $u > A$, $v > A$. For $\ell(\phi'(S_{t+j}))$ with $t + j \leq \tau_{Ax}$ we have then the same inequality as in (iii), thus, the contribution of such terms can be estimated as in the case (i), but with the exponent $\alpha - \delta$. The remaining sum is bounded from above by

$$\vec{\nu}(A)\mathbb{E}\sigma_{1-q(x)}$$

analogously to the case (ii).   $\square$

The integrability condition (12) can be ensured by means of the following lemma found in [3], page 28.

LEMMA 3.4.   *For each decreasing positive function $g$ on $[0,\infty[$ and each subordinator $(S_t, t \geq 0)$, the following properties are equivalent:*

$$\int_0^\infty g(t)\,dt < \infty \quad \Longleftrightarrow \quad \int_0^\infty g(S_t)\,dt < \infty \qquad a.s.$$

$$\Longleftrightarrow \quad \mathbb{E}\int_0^\infty g(S_t)\,dt < \infty.$$

REMARKS. The function $\ell(1/y) = |\log y|^\rho, \rho > 0$, is decreasing, thus, part (i) of the theorem applies. For $\phi(y) = 1 - (y + 1)^{-\beta}$, the integrability condition is fulfilled if $\alpha(\beta + 1) > 1$, thus, selecting $b$ in the range $1 - \alpha < b < \alpha\beta$, part (iv) applies for any $\ell$, with $q(y) = y^b$. Part (ii) is useful for $\ell(1/y) = |\log y|^{-\rho}$, $\rho > 0$, in which case we can take $r(y) = y^b, 0 < b < 1$, to meet (iib).

Note that, for $q(y) = y^b, 0 < b < 1$, the condition on $\phi^\leftarrow$ can be reformulated in terms of $\phi$, for example, $\phi^\leftarrow(1 - x) = o(x^{-\alpha})$ is equivalent to $1 - \phi(y) = o(y^{-1/\alpha})$.

**4. Strong laws.** Strong laws for $K_n$ and the like for composition derived from a transformed subordinator follow by combining the results in the two previous sections. Introduce

$$(15) \qquad L(t) = \int_0^t (\phi'(S_u))^\alpha \, du, \qquad L = L(\infty).$$

Recall that notation $K_n(t)$ or $K_{n,r}(t)$ refers to the parts of the partial composition produced by the range of $(\phi(S_u), u \in [0, t])$.

THEOREM 4.1. *For $0 < \alpha < 1$, the regular variation assumption* (1) *implies*

$$\frac{K_n(t)}{\Gamma(1 - \alpha)n^\alpha \ell(n)} \to L(t), \qquad \frac{K_{n,r}(t)}{\Gamma(1 - \alpha)n^\alpha \ell(n)} \to (-1)^{r-1} \binom{\alpha}{r} L(t),$$

*with probability 1, as $n \to \infty$. And if $\phi$ satisfies also the conditions of Theorem* 3.3, *we have*

$$\frac{K_n}{\Gamma(1 - \alpha)n^\alpha \ell(n)} \to L, \qquad \frac{K_{n,r}}{\Gamma(1 - \alpha)n^\alpha \ell(n)} \to (-1)^{r-1} \binom{\alpha}{r} L.$$

*For $\alpha = 1$, the analogous results are read from Theorem* 2.1.

*Asymptotics in the drift case.* We sketch the extension to the case of subordinator with positive drift $\mathtt{d}$. In this case the length of composition satisfies $K_n \sim K_{n,1} \sim n\lambda(\bar{\mathcal{R}})$, where $\lambda$ denotes the Lebesgue measure, thus, the asymptotics follow from the next lemma which generalizes [16], Corollary 5.8.

LEMMA 4.2. *Let $(S_t)$ be a subordinator with drift $\mathtt{d} > 0$, and $\phi : [0, \infty] \to [0, 1]$ be an absolutely continuous strictly increasing function. Then*

$$\lambda(\phi(\mathcal{R})) = \mathtt{d} \int_0^\infty \phi'(S_t) \, dt.$$

PROOF.    Because $(S_t)$ is almost everywhere differentiable with derivative $\mathtt{d}$, by the change of variable $S_t = y$,

$$\mathtt{d} \int_0^\infty \phi'(S_t) \, dt = \int_{\mathcal{R}} \phi'(y) \, dy = \lambda(\phi(\mathcal{R})).  \qquad \square$$

In the drift case it is sensible to distinguish between *genuine* singleton parts which are caused when a sample point hits $\widetilde{\mathcal{R}}$, and the other *occasional* singletons induced by open gaps comprising $\widetilde{\mathcal{R}}^c$. (For fixed $n$, the composition $\mathcal{C}_n$ does not allow one to make this distinction.) Denoting $K_{n,1-}$ and $K_{n,1+}$ the counts of genuine and occasional singletons, we have $K_n = K_{n,1-} + K_{n,1+} + \sum_{r=2}^\infty K_{n,r}$, and $K_n \sim K_{n,1-}$. The asymptotic behavior of variables $K_{n,1+}$ and $K_{n,r}$ for $r > 1$ is then still as in Theorem 2.1, as follows by noting that the gap-counting Theorems 3.1 and 3.3 and Corollary 3.2 are also valid for subordinators with drift.

A curious phenomenon occurs for $\alpha = 1$: normalizing the Lévy data $(\nu, \mathtt{d})$ so that $\mathtt{d} = 1$, we have all three variables $K_n/n$, $K_{n,1-}/n$ and $K_{n,1+}/(n\ell^*(n))$ approaching the same limit.

## 5. Regenerative compositions.
We shall specialize the results of the previous section to the regenerative compositions appearing when $\widetilde{\mathcal{R}}$ is the range of the multiplicative subordinator $\widetilde{S}_t = 1 - \exp(-S_t)$. In this case there is a simple connection between $\widetilde{S}_t$ and $L(t)$ and a nice formula for the moments of $L$.

With each multiplicative subordinator, we associate the *area process*

$$A(t) = \int_0^t (1 - \widetilde{S}_u) \, du$$

and its terminal value $A = A(\infty)$ obtained by taking the infinite integration bound. In terms of the (additive) subordinator,

$$A = \int_0^\infty \exp(-S_u) \, du = \int_0^\infty (1 - \widetilde{S}_u) \, du$$

is the widely studied *exponential functional*, see [4, 5, 7, 28, 29].

Let $\widetilde{\nu}$ be the measure on $[0,1]$ obtained by the exponential transform $\phi(y) = 1 - e^{-y}$ from the measure $\nu$. The Laplace exponent is thus given by

$$\Phi(s) = \int_0^\infty (1 - e^{-sy}) \nu(dy) = \int_0^1 (1 - (1-x)^s) \widetilde{\nu}(dx).$$

Because $\phi'(0) = 1$, the assumption (1) implies that the tail of the transformed measure satisfies $\widetilde{\nu}[x, 1] \sim x^{-\alpha} \ell(1/x)$ for $x \downarrow 0$. For arbitrary $\alpha > 0$, the process $\widetilde{S}_t^{(\alpha)} := 1 - (1 - \widetilde{S}_t)^\alpha$ is itself a multiplicative subordinator with Lévy measure $\widetilde{\nu}^{(\alpha)}$ related to $\widetilde{\nu}$ by $\widetilde{\nu}^{(\alpha)}[x,1] = \widetilde{\nu}[1-(1-x)^{1/\alpha},1]$. The relation between the corresponding Laplace exponents is $\Phi^{(\alpha)}(s) = \Phi(\alpha s)$. That is to say,

$$\widetilde{S}_t^{(\alpha)} = 1 - \exp(-\alpha S_t)$$

is the multiplicative counterpart of the scaled subordinator $(\alpha S_t)$. The area process for $(\widetilde{S}_t^{(\alpha)})$ is

$$L^{(\alpha)}(t) := \int_0^t (1 - \widetilde{S}_u^{(\alpha)})\, du = \int_0^t (1 - \widetilde{S}_u)^\alpha\, du = \int_0^t \exp(-\alpha S_u)\, du$$

and we define $L^{(\alpha)} = L^{(\alpha)}(\infty)$ to be the $A$-functional for $(S_t^{(\alpha)})$.

For the scaling function $\psi(x) = x^{-\alpha}\ell(1/x)$, we have the following:

THEOREM 5.1. *Suppose the Lévy measure fulfills* (1), *then, for $x \downarrow 0$, the jump counts of the multiplicative subordinator $\widetilde{S}_t = 1 - \exp(-S_t)$ satisfy*

$$\widetilde{N}_x(t)/\psi(x) \to L^{(\alpha)}(t), \qquad \widetilde{N}_x/\psi(x) \to L^{(\alpha)},$$

*for $x \downarrow 0$ almost surely and in the mean.*

PROOF. We have $\phi'(y) = e^{-y}$. Then part (iv) of Theorem 3.3 applies, because the required integrability holds for any arbitrary positive power and the second condition is fulfilled with $q(x) = x$. □

The application of Theorem 2.1 results in the following:

COROLLARY 5.2. *If the Lévy measure fulfills* (1), *then, for $0 < \alpha < 1$, the composition induced by the multiplicatively regenerative set $\widetilde{\mathcal{R}}$ satisfies, for $n \to \infty$,*

$$K_n/(n^\alpha \ell(n)) \to \Gamma(1-\alpha)L^{(\alpha)}$$

*almost surely and in the mean. And for $\alpha = 1$,*

$$K_n/(n\ell^*(n)) \to L^{(1)} = A$$

*almost surely and in the mean, with $\ell^*$ as in* (3).

Generalizations to $K_n(t)$ and $K_{n,r}(t)$ follow in the same way.

The distribution of $L^{(\alpha)}$ admits some exponential moments, hence, it is uniquely determined by its integer moments. They are given by the following formula which was recorded for general Lévy processes in [7], though it can be traced back in special cases to much earlier literature (see, e.g., [10], page 283):

$$(16) \qquad \mathbb{E}(L^{(\alpha)})^k = \frac{k!}{\prod_{j=1}^k \Phi(\alpha j)}, \qquad k = 1, 2, \ldots.$$

Example. Consider the two-parameter family of regenerative composition structures, as in [16], with

$$\tilde{\nu}[x,1] = x^{-\alpha}(1-x)^{\theta}, \qquad 0 < \alpha < 1, \theta > 0.$$

In this case we have

$$\Phi(s) = \frac{s\Gamma(1-\alpha)\Gamma(s+\theta)}{\Gamma(s+1-\alpha+\theta)},$$

thus,

$$\mathbb{E}(L^{(\alpha)})^k = \frac{\Gamma(\theta+1)(\alpha+\theta)(2\alpha+\theta)\cdots((k-1)\alpha+\theta)}{\Gamma(k\alpha+\theta)\alpha^k}$$

in agreement with formula (192) from [25]. This specializes for $\theta = 0$ to the integer moments of the Mittag–Leffler distribution with parameter $\alpha$.

Now we shall give a new proof of (16) in the case of subordinators. The method we use here is not applicable to the exponential functionals of more general Lévy processes, as considered in [5, 7], but it is much more elementary and apparently more natural in the context of multiplicatively regenerative sets. By the above discussion, it is sufficient to prove the formula for $\alpha = 1$, that is, for the area functional of a multiplicative subordinator. Letting $m_k = \mathbb{E}A^k$, we wish to show that

$$(17) \qquad m_k = \frac{k!}{\prod_{j=1}^k \Phi(j)}, \qquad k = 0, 1, \dots.$$

*Finite Lévy measure, no drift.* Suppose first that $\mathtt{d} = 0$ and $\tilde{\nu}$ is a probability measure. Let $X_j$ be a sample from $\tilde{\nu}$. Then $(\widetilde{S}_t)$ is a step function whose range is a stick-breaking sequence $X_j(1-X_1)\cdots(1-X_{j-1}), j = 1, 2, \dots,$ complemented by 0 and 1. The jumps of $(\widetilde{S}_t)$ occur at the epochs of an independent homogeneous Poisson process. Therefore, $A$ is representable as a random series

$$A = E_1 X_1 + (E_1 + E_2)X_2(1-X_1) + (E_1 + E_2 + E_3)X_3(1-X_1)(1-X_2) + \cdots,$$

where $E_j$ are jointly independent exponential random variables with mean 1, also independent from the $X_j$'s. The series is finite only if $\tilde{\nu}\{1\} > 0$. We can also re-arrange terms and write $A$ in the form

$$A = E_1 + E_2(1-X_1) + E_3(1-X_1)(1-X_2) + \cdots,$$

from which we deduce

$$(18) \qquad A \stackrel{d}{=} E + (1-X)A',$$

where $A'$ is a replica of $A$, independent of $(E, X) \overset{d}{=} (E_1, X_1)$. By virtue of the formula

$$\mathbb{E}(1 - X)^k = 1 - \Phi(k),$$

the expectation is computed as $\mathbb{E}A = 1/\Phi(1)$, in accord with the $k = 1$ case of (17). Furthermore, the identity

$$A^k \overset{d}{=} \sum_{j=0}^{k} \binom{k}{j} E^{k-j} (1 - X)^j A^j$$

implies that the moments $(m_k)$ of $A$ exist and satisfy the recursion

$$m_k \Phi(k) = \sum_{j=0}^{k-1} \frac{k!}{j!} (1 - \Phi(j)) m_j.$$

To solve the recursion, split out the last term and substitute the same identity but with $k - 1$ to arrive at

$$m_k \Phi(k) = k(1 - \Phi(k-1)) m_{k-1} + k m_{k-1} \Phi(k-1),$$

which is the same as the simple multiplicative recursion

$$m_k \Phi(k) = k m_{k-1},$$

whose unique solution with the initial value $m_0 = 1$ is that given by (17).

Replacing the probability measure $\widetilde{\nu}$ by its positive multiple, say, $\widetilde{\nu}_c = c\widetilde{\nu}$, we can write a series representation for the corresponding functional $A_c$ exactly as above, but with $E_j/c$ instead of $E_j$. Since $\mathbb{E}(E/c)^k = k!/c^k$, the same computation yields the additional factor $c^k$ in the denominator. This agrees with (17) because the new Laplace exponent is the multiple $c\Phi$.

*Finite Lévy measure, positive drift.* The moments formula for a subordinator with drift can be proved analytically using approximation by drift-free subordinators, but it is more instructive to inquire into this case separately. Assuming $\widetilde{\nu}[0, 1] = 1$ and $\mathtt{d} > 0$, the subordinator coincides with function $1 - e^{-\mathtt{d}t}$ for $t < E_1$, and the jump at time $E_1$ is $e^{-\mathtt{d}E_1} X_1$. The first-jump decomposition is

$$A \overset{d}{=} e^{-\mathtt{d}E}(1 - X) A' + (1 - e^{-\mathtt{d}E})/\mathtt{d},$$

with the same convention as in (18). The recursion for moments is obtained by using

$$\mathbb{E}(1 - X)^k = 1 - \Phi(k) + k\mathtt{d}$$

and exploiting the fact that $e^{-\mathtt{d}E}$ is distributed according to $\mathrm{Beta}(1/\mathtt{d},1)$, we obtain

$$m_k \frac{\Phi(k)}{k\mathtt{d}+1} = \sum_{j=0}^{k-1} \frac{k!}{j!} \frac{(1-\Phi(j)+j\mathtt{d})}{\mathtt{d}^{k-j+1}} \frac{\Gamma(j+1/\mathtt{d})}{\Gamma(k+1+1/\mathtt{d})} m_j.$$

The solution is again (17), as justified by the same inductive argument. Repeating the above scaling argument, we see that the formula also holds in the case $\mathtt{d}>0$ and arbitrary finite $\widetilde{\nu}$.

*General subordinator.* Given arbitrary Lévy data $(\mathtt{d},\widetilde{\nu})$, consider the family of subordinators with parameters $(\mathtt{d},\widetilde{\nu}_\varepsilon)$, where $\widetilde{\nu}_\varepsilon$ is a truncated measure that coincides with $\widetilde{\nu}$ outside $[0,\varepsilon]$ and is zero within $[0,\varepsilon]$. Using a version of the well-known recipe, we can construct the corresponding multiplicative subordinators $\widetilde{S}_t$ and $\widetilde{S}_{\varepsilon,t}$ using the same Poisson point process in the strip $[0,\infty[\times[0,1]$ with intensity measure Lebesgue $\times\,\widetilde{\nu}$, so that

$$\widetilde{S}_t = 1 - e^{-\mathtt{d}t} \prod_{\tau_j \leq t} (1-\Delta_j),$$

where the product is over the atoms $(\tau_j,\Delta_j)$, and $\widetilde{S}_{\varepsilon,t}$ has a similar representation, with the only distinction that the factors corresponding to the atoms with $\Delta_j \leq \varepsilon$ do not enter into the product. By construction, $\widetilde{S}_{\varepsilon,t} \uparrow \widetilde{S}_t$ as $\varepsilon \downarrow 0$ and the convergence is uniform in $t$. Thus, by monotone convergence, we have for the corresponding integrals $A_\varepsilon(t_1,t_2) \downarrow A(t_1,t_2)$ almost surely and with all moments, for all $0 \leq t_1 < t_2 \leq \infty$. But all measures $\widetilde{\nu}_\varepsilon$ are finite, thus, as we have shown, the moments formula is true, therefore, the formula also holds for $(\mathtt{d},\widetilde{\nu})$ since the corresponding Laplace exponents satisfy $\Phi_\varepsilon \uparrow \Phi$.

REMARK. For $A'$ a copy of $A$ independent of $(\widetilde{S}_t,A(t))$ (fixed $t$), we have

$$A(t,\infty) \overset{d}{=} (1-\widetilde{S}_t)A', \qquad A \overset{d}{=} A(t) + (1-\widetilde{S}_t)A'.$$

This leads to

$$\mathbb{E}(A(t,\infty))^k = \frac{k!}{\prod_{j=1}^k \Phi(j)} \exp(-t\Phi(k)), \qquad \mathbb{E}A(t) = \frac{1-\exp(-t\Phi(1))}{\Phi(1)}.$$

Higher moments of $A(t)$ are not immediate because of dependence between $S_t$ and $A(t)$. See [7] for the formulas with $t$ replaced by an independent exponential variable.

**6. Poissonized compositions.** A closely related type of structure appears when the uniform sample of fixed size $n$ is replaced by a Poisson point process of rate $\rho$. A composition of random weight $n(\rho)$ appears by separating the Poisson points into blocks by means of a random closed set $\widehat{\mathcal{R}} \subset [0, 1]$. We shall denote this poissonized composition $\widehat{\mathcal{C}}_\rho$, and provide with "ˆ" all quantities related to it. The relation to the fixed-$n$ composition is therefore $\widehat{\mathcal{C}}_\rho \overset{d}{=} \mathcal{C}_n$, conditionally, given $n(\rho) = n$.

Poissonization is useful for two reasons. Generally speaking, it is a powerful technique for asymptotic considerations in combinatorial problems, allowing one to explicitly exploit the spatial independence where otherwise only a kind of asymptotic independence is available (see, e.g., [27] for overview). On the other hand, poissonization yields a family of compositions $(\mathcal{C}_\rho, \rho > 0)$ which satisfies a consistency condition analogous to the defining property of partition or composition structures [11, 16, 25]. Explicitly, for any $\rho > 0$ and $x \in [0, 1]$, probability distributions of the following compositions coincide: (i) a composition with rate parameter $\rho(1 - x)$ and (ii) a *thinned* composition which appears when the atoms making up a sample with rate $\rho$ are deleted independently with probability $x$.

In the sequel we shall only consider the case when $\widetilde{\mathcal{R}}$ is the range of a multiplicative subordinator $\widetilde{S}_t = 1 - \exp(-S_t)$, in which case the distribution of (i) also coincides with the distribution of (iii), a tail composition of the composition of rate $\rho$ which appears to the right of $x$, conditionally, given $x \in \widetilde{\mathcal{R}}$. The last equivalence is analogous to the regenerative property of the fixed-$n$ compositions [16]. The same composition of random integer $n(\rho)$ appears when the range of additive subordinator $(S_t)$ is used to separate into blocks atoms of inhomogeneous Poisson process on $[0, \infty)$ with exponential intensity measure $\rho e^{-x} dx, x \in [0, \infty]$. We denote by $\widehat{K}_\rho$ the length of the poissonized composition, and $\widehat{K}_\rho(t), \widehat{K}_{\rho,r}(t)$ stand for the number of parts of the partial composition produced by the range of multiplicative subordinator up to time $t$.

We proceed with the convergence results which recover and complement the results in the previous sections. The equivalence of strong laws for $(\widehat{\mathcal{C}}_\rho)$ and $(\mathcal{C}_n)$ is quite obvious, and for quantities like moments, there is a well-developed analytical technique of poissonization/depoissonization [19], though we shall use more elementary arguments.

6.1. *Recursions.* Let $\mathcal{F}_t$ be the $\sigma$-algebra generated by the subordinator $(S_u, u \in [0, t])$ and by the Poisson configuration on $[0, S_t]$. By the independence property of the Poisson process and the regenerative property of $\widetilde{\mathcal{R}}$, the tail composition induced by $\widetilde{\mathcal{R}} \cap ]\widetilde{S}_t, 1]$ is independent of $\mathcal{F}_t$. This observation is a source of recursions related to the poissonized composition.

Let $p_j(\rho), j = 0, 1, \ldots$, be the distribution of $\widehat{K}_\rho$. Each $p_j$ may be extended to an entire function of the complex variable, with initial value $p_j(0) = 0$ [with the only exception $p_0(\rho) = e^{-\rho}$]. Introduce the factorial moments $f^{(m)}(\rho) = \mathbb{E}\widehat{K}_\rho(\widehat{K}_\rho - 1)\cdots(\widehat{K}_\rho - m + 1)$, with $m = 0, 1, \ldots, f^{(0)}(\rho) = 1$.

LEMMA 6.1. *The following integral recursions hold. For $j = 1, 2, \ldots$,*

(19)
$$\int_0^1 (p_j(\rho) - e^{-\rho x} p_j(\rho(1-x))\widetilde{\nu}(dx))$$
$$= \int_0^1 (1 - e^{-\rho x})p_{j-1}(\rho(1-x))\widetilde{\nu}(dx),$$

*and for $m = 1, 2, \ldots$,*

(20)
$$\int_0^1 (f^{(m)}(\rho) - f^{(m)}(\rho(1-x)))\widetilde{\nu}(dx)$$
$$= m \int_0^1 (1 - e^{-\rho x})f^{(m-1)}(\rho(1-x))\widetilde{\nu}(dx).$$

PROOF. Each $p_j(\rho)$ may be written as a generating function whose coefficients are rational functions in the variables $\Phi(n), n \geq 0$, for example, the probability of one-part composition is

$$p_1(\rho) = e^{-\rho} \sum_{n=1}^\infty \frac{\rho^n}{n!} \frac{\Phi(n:n)}{\Phi(n)},$$

where

$$\Phi(n:m) = \binom{n}{m} \int_0^1 x^m (1-x)^n \widetilde{\nu}(dx), \qquad 1 \leq m \leq n.$$

Thus, the statement is of purely algebraic nature and can be translated as a series of polynomial identities in these variables. Thus, it is sufficient to consider the "stick-breaking case" of finite Lévy measure, normalized to a probability measure, when the recursion is proved by conditioning on the first break $X = x$ with distribution $\widetilde{\nu}$. Indeed, there are $j$ blocks when either $[0, x]$ contains a Poisson atom and then $[x, 1]$ generates $j - 1$ blocks [with probability $p_{j-1}(\rho(1-x))$], or $[0, x]$ is empty and $[x, 1]$ generates $j$ blocks. This gives

$$p_j(\rho) = \int_0^1 [e^{-\rho x} p_j(\rho(1-x)) + (1 - e^{-\rho x})p_{j-1}(\rho(1-x))]\widetilde{\nu}(dx).$$

To keep this formula right for arbitrary finite $\widetilde{\nu}$, we should put $p_j(\rho)$ into the integral, then the formula becomes homogeneous in the $\Phi(n)$'s and holds in general, for algebraic reasons.

To prove (20), start with the definition

$$f^{(m)}(\rho) = \sum_{j=0}^{\infty} p_j(\rho) j(j-1)\cdots(j-m+1),$$

then multiply both parts in (19) by $j(j-1)\cdots(j-m+1)$ and sum over $j$ using the identity

$$j(j-1)\cdots(j-m+1) = (j-1)\cdots(j-m) + m(j-1)\cdots(j-m+1). \quad \square$$

Manipulation with power series allows, in principle, computing the distribution of $K_n$ with all moments, for fixed-$n$ compositions. Let us demonstrate this on the expectation $f^{(1)}(\rho) = \mathbb{E}\widehat{K}_\rho$. For $f_n^{(1)} = \mathbb{E}K_n$, we have the poissonization identity

$$f^{(1)}(\rho) = \sum_{n=0}^{\infty} e^{-\rho} \frac{\rho^n}{n!} f_n^{(1)}.$$

Substituting into (20) and integrating, we obtain a relation between generating functions

$$e^{-\rho} \sum_{n=0}^{\infty} \frac{\rho^n}{n!} f_n^{(1)} \Phi(n) = \sum_{n=1}^{\infty} (-1)^n \frac{\rho^n}{n!} \Phi(n:n).$$

Multiplying by $e^\rho$ and extracting the coefficients, we get

$$f_n^{(1)} = \sum_{j=1}^{n} (-1)^{j+1} \binom{n}{j} \frac{\Phi(j:j)}{\Phi(j)},$$

which is a familiar expression for $\mathbb{E}K_n$, see [16].

6.2. *Asymptotics.* The convergence $\widehat{K}_\rho/(\Gamma(1-\alpha)\rho^\alpha \ell(\rho)) \to L^{(\alpha)}$ a.s. follows exactly as in Sections 4 and 5 for $0 < \alpha < 1$ (and with a proper scaling, also for $\alpha = 1$), in the footprints of Karlin [20], where the Poisson model was treated in parallel with the fixed-$n$ case. In this section we show that the recursion (20) implies the moments formulae analogous to (16). This can be regarded as a proof that the convergence holds together with all moments, and also as yet another derivation of the moments formula (16).

Introducing the poissonized Laplace exponent (not to be confused with $\Phi$ written in terms of $\nu$)

$$\widehat{\Phi}(\rho) = \int_0^1 (1 - e^{-\rho x}) \widetilde{\nu}(dx),$$

we have as $\rho \to \infty$

$$\widehat{\Phi}(\rho) \sim \begin{cases} \Gamma(1-\alpha)\rho^\alpha \ell(\rho), & \text{for } 0 < \alpha < 1, \\ \rho \ell^*(\rho), & \text{for } \alpha = 1, \end{cases}$$

(see the Appendix).

PROPOSITION 6.2. *The factorial moments satisfy*

$$f^{(m)}(\rho) \sim c^{(m)}(\rho^\alpha \ell(\rho))^m,$$

*with $c^{(m)}$ given by*

(21)
$$c^{(m)} = \prod_{j=1}^{m} \frac{j\Gamma(1-\alpha)}{\Phi(\alpha j)}$$

*for $0 < \alpha < 1$. For $\alpha = 1$, the factors $\Gamma(1-\alpha)$ should be omitted and $\ell$ replaced by $\ell^*$.*

PROOF. We concentrate on the case $0 < \alpha < 1$, leaving the case $\alpha = 1$ to the reader. Trivially, $f^{(0)}(\rho) = 1$. Suppose by induction that the asymptotics holds for some $m - 1$. Then setting $b = \Gamma(1-\alpha)$ and $g(\rho) = (\rho^\alpha \ell(\rho))^m$, we have

(22)
$$\int_0^1 (1 - e^{-\rho x}) f^{(m-1)}(\rho(1-x)) \widetilde{\nu}(dx) \sim bc^{(m-1)} g(\rho).$$

This is shown by splitting the integral at $\varepsilon$ and replacing the integrand for $x \in ]0, \varepsilon]$ by its asymptotics. To justify the induction step, fix $\varepsilon$ and suppose there exists arbitrarily large $\rho$ such that

$$f^{(m)}(\rho) > (1 + \varepsilon) c^{(m)} g(\rho)$$

(we wish to lead this assumption to a contradiction). Then, perhaps selecting $\varepsilon$ smaller, for any fixed constant $C$, there exists arbitrarily large $\rho$ such that

$$f^{(m)}(\rho) > (1 + \varepsilon) c^{(m)} g(\rho) + C$$

[just because $f^{(m)}(\rho) \to \infty$]. Up to the end of this paragraph, $\rho = \rho(C)$ will be the minimal $\rho$ for which the inequality holds. Note that $\rho(C) \to \infty$ as $C \to \infty$. Thus, we have

$$f^{(m)}(\rho) = (1 + \varepsilon) c^{(m)} g(\rho) + C,$$

$$f^{(m)}(\rho(1-x)) < (1 + \varepsilon) c^{(m)} g(\rho)(1-x)^{m\alpha} + C, \qquad x \in ]0, 1],$$

and substituting this into (20), we see that the left-hand side is estimated from below by

$$(1 + \varepsilon) c^{(m)} g(\rho) \int_0^1 (1 - (1-x)^{\alpha m}) \widetilde{\nu}(dx) = (1 + \varepsilon) c^{(m)} g(\rho) \Phi(\alpha m)$$

$$= (1 + \varepsilon) c^{(m-1)} b g(\rho),$$

where $\rho \to \infty$ and we used monotonicity. This disagrees with the right-hand side of (20) given by (22), giving the required contradiction. Thus,

$$\limsup \frac{f^{(m)}(\rho)}{c^{(m)} g(\rho)} < 1 + \varepsilon.$$

A symmetric argument shows that

$$\liminf \frac{f^{(m)}(\rho)}{c^{(m)}g(\rho)} > 1 - \varepsilon.$$

Letting $\varepsilon \to 0$ ends the proof. □

Depoissonization follows rather easily. Recall that the collection of atoms of the Poisson process with rate $\rho$ can be identified with a uniform sample $\{u_1, \ldots, u_{n(\rho)}\}$, with $n(\rho)$ distributed according to Poisson $(\rho)$. By the obvious monotonicity of $K_n$ we have

$$K_n \mathbf{1}(n(\rho(1-\varepsilon)) < n) \geq \widehat{K}_{\rho(1-\varepsilon)} \mathbf{1}(n(\rho(1-\varepsilon)) < n),$$

$$K_n \mathbf{1}(n(\rho(1+\varepsilon)) > n) \leq \widehat{K}_{\rho(1+\varepsilon)} \mathbf{1}(n(\rho(1+\varepsilon)) > n).$$

Selecting $n = \rho$ and letting $\rho \to \infty$, the elementary large deviation bounds for the probability $\mathbb{P}(\rho(1-\varepsilon) < n(\rho) < \rho(1+\varepsilon))$ imply that, for $n \to \infty$,

$$\widehat{K}_{n(1-\varepsilon)} \ll K_n \ll \widehat{K}_{n(1+\varepsilon)}.$$

Letting $\varepsilon \downarrow 0$ and using Proposition 6.2, we see that $K_n \sim \widehat{K}_n$ almost surely and with all moments.

Observing that the computation of the constants (21) is equivalent to the formula (16) [thus, we have yet another proof for (16)], and recalling Theorem 2.1, we summarize the above discussion in the following theorem.

THEOREM 6.3. *The almost-sure convergence*

$$\widehat{K}_\rho/(\rho^\alpha \ell(\rho)) \to \Gamma(1-\alpha)L^{(\alpha)}, \qquad \rho \to \infty,$$

$$K_n/(n^\alpha \ell(n)) \to \Gamma(1-\alpha)L^{(\alpha)}, \qquad n \to \infty,$$

*holds for both Poisson and fixed-n compositions together with the convergence of all integer moments for $0 < \alpha < 1$, and with a proper scaling also for $\alpha = 1$.*

The convergence of moments of $K_{n,r}, K_n(t), K_{n,r}(t), \widehat{K}_{\rho,r}, \widehat{K}_\rho(t), \widehat{K}_{\rho,r}(t)$ (with obvious definition of the last two random variables) follows from the theorem by dominated convergence.

6.3. *A martingale approach.* Extending the discussion in the previous section, consider $\widehat{K}_\rho(t)$, which is the number of blocks of a poissonized composition produced by the subordinator up to time $t$. We can view $(\widehat{K}_\rho(t), t \geq 0)$ as either an increasing process with unit jumps or a point process of those jump-times of $(\widetilde{S}_t)$ which have jump intervals covering some sample points. The compensator for $\widehat{K}_\rho(t)$ is

$$C_t = \int_0^t \widehat{\Phi}(\rho(1-\widetilde{S}_u)) \, du.$$

By observing that $\widehat{\Phi}(\rho(1-x))$ is the probability that a gap with leftpoint $x$ is hit by a Poisson atom, the formula can be first argued in the renewal case. The general case follows by extrapolation from the case of finite Lévy measure. This readily implies the following:

LEMMA 6.4.  *For each $\rho > 0$, the process*

$$M_t := \widehat{K}_\rho(t) - C_t, \qquad t \in [0, \infty],$$

*is a square-integrable martingale with unit jumps and quadratic predictable characteristics $\langle M \rangle_t = C_t$. Furthermore,*

$$\mathbb{E}M_t^2 = \mathbb{E}(\widehat{K}_\rho(t) - C_t)^2 = \mathbb{E}C_t.$$

PROOF.  The squared jump magnitudes of $M_t$ are 1. This implies that the submartingale $M_t^2$ has the same compensator as $\widehat{K}_\rho(t)$, that is, $\langle M \rangle_t = C_t$, as in [22], Section 6.2.  □

The lemma opens yet another approach to the convergence results, for which we give below the $\mathcal{L}^2$-version. Note that the scaling by $\Phi(\rho)$ is asymptotically the same as the scaling by $\widehat{\Phi}(\rho)$ (see Appendix), which is asymptotic to $\Gamma(1-\alpha)\ell(\rho)\rho^\alpha$ for $0 < \alpha < 1$ and to $\ell^*(\rho)\rho$ for $\alpha = 1$.

THEOREM 6.5.  *Under the regular variation assumption* (1), *as $\rho \to \infty$,*

$$(23) \qquad \frac{\widehat{K}_\rho}{\Phi(\rho)} \to L^{(\alpha)}$$

*almost surely and in $\mathcal{L}^2$. An analogous result is valid for $\widehat{K}_\rho(t)$ for each $t > 0$.*

PROOF.  We wish to establish the convergence (23) in $\mathcal{L}^2$. Use Lemma 6.4 to obtain

$$(24) \qquad \mathbb{E}\left(\frac{\widehat{K}_\rho}{\widehat{\Phi}(\rho)} - \int_0^\infty \frac{\widehat{\Phi}(\rho(1-\widetilde{S}_t))}{\widehat{\Phi}(\rho)}\, dt\right)^2 = \frac{1}{\widehat{\Phi}(\rho)}\mathbb{E}\int_0^\infty \frac{\widehat{\Phi}(\rho(1-\widetilde{S}_t))}{\widehat{\Phi}(\rho)}\, dt.$$

Also observe that

$$(25) \qquad \frac{\widehat{\Phi}(\rho(1-\widetilde{S}_t))}{\widehat{\Phi}(\rho)} \to (1-\widetilde{S}_t)^\alpha \qquad \text{as } \rho \to \infty$$

almost surely for each fixed $t$. Thus, (23) would follow by dominated convergence once we could bound

$$\int_0^\infty \frac{\widehat{\Phi}(\rho(1-\widetilde{S}_t))}{\widehat{\Phi}(\rho)}\, dt$$

from above by a square-integrable random variable. To this end, write $\widehat{\Phi}(\rho) = \rho^\alpha \ell_0(\rho)$ with slowly varying $\ell_0$ [so $\ell_0(\rho) \sim \ell(\rho)$], then, if the Potter's bound were valid for $\ell_0$ on $]0, \infty[$, we could estimate

$$\widehat{\Phi}(\rho(1 - \widetilde{S}_t))/\widehat{\Phi}(\rho) < C(1 - \widetilde{S}_t)^{\alpha - \delta}$$

with some small $\delta > 0$ and $C > 0$, whence

$$\int_0^\infty (\widehat{\Phi}(\rho(1 - \widetilde{S}_t))/\widehat{\Phi}(\rho)) \, dt < C \int_0^\infty (1 - \widetilde{S}_t)^{\alpha - \delta} \, dt \in \mathcal{L}^2.$$

To make this argument precise, we fix some sufficiently large constant $c = X(C, \delta)$ required in Theorem A.3(i), and then split the integral at the first passage time $\sigma_{1 - c/\rho}$ of $(\widetilde{S}_t)$ over the level $1 - c/\rho$. The tail integral

$$\int_{\sigma_{1 - c/\rho}}^\infty \widehat{\Phi}(\rho(1 - \widetilde{S}_t)) \, dt$$

is bounded from above by a Poisson random variable with mean $c$, which is obviously square-integrable, for each $c > 0$. Thus, it is sufficient to exploit Potter's bound for $\ell_0$ on $[c, \infty[$.  $\square$

## APPENDIX

**A.1. Abel–Tauberian theorems.** An exposition of Abel–Tauberian theorems for the Laplace transform is given in [9], Section XIII.5. Bingham, Goldie and Teugels ([6], Section 1.7 and Chapter 4) give a fuller account, also for more general integral transforms, including the Mellin transform.

We establish next some elementary connections between the integral transforms in a form suitable for applications to subordinators. Consider the two Laplace exponents (also called Bernstein functions)

$$\Phi(s) := \int_0^\infty (1 - e^{-sy}) \nu(dy) = \int_0^1 (1 - (1 - x)^s) \widetilde{\nu}(dx)$$

and

$$\widehat{\Phi}(s) := \int_0^1 (1 - e^{-sx}) \widetilde{\nu}(dx) = \int_0^\infty (1 - e^{-s(1 - e^{-y})}) \nu(dy),$$

where the measure $\widetilde{\nu}(dx)$ on $]0, 1]$ is the image of $\nu(dy)$ on $]0, \infty]$ via $x = 1 - e^{-y}$, and it is assumed that $\Phi(s) < \infty$ for some (and, hence, all) $s > 0$. The function $\widehat{\Phi}$ is the poissonization of $\Phi$, that is,

$$\widehat{\Phi}(s) = \sum_{n=0}^\infty e^{-s} \frac{s^n}{n!} \Phi(n).$$

LEMMA A.1.   *Whatever the Laplace exponent* $\Phi$,

$$\lim_{s\to\infty}(\Phi(s) - \widehat{\Phi}(s)) = 0.$$

PROOF.   With a hint from [8], page 1257, we have, for $0 \le x \le 1$ and $s \ge 0$,

$$0 \le e^{-sx} - (1-x)^s \le sx^2 e^{-sx} \le e^{-1}x.$$

Using this we get

$$0 \le \Phi(s) - \widehat{\Phi}(s) = \int_0^1 (e^{-sx} - (1-x)^s)\tilde{\nu}(dx) < e^{-1}\int_0^1 x\tilde{\nu}(dx) = e^{-1}\Phi(1),$$

so the claim follows by dominated convergence.   □

COROLLARY A.2.   *Abel–Tauberian relations (as* $s \to \infty$*) between the tail of the measure* $\nu[1/s, \infty]$ *and the Laplace exponent are the same for* $\Phi(s)$ *and* $\widehat{\Phi}(s)$.

Estimates of the difference can be given under the assumption of regular variation. The Laplacian case of monotone density, relating asymptotics of $\tilde{\nu}$ and $\widehat{\Phi}$, is covered by a combination of Theorems 3 and 4 in [9], Section 13.5.

## A.2. Potter's theorem.

THEOREM A.3 ([6], Theorem 1.5.6).   *Let* $\ell$ *be a function of slow variation at infinity.*

(i) *For arbitrarily chosen constants* $A > 1, \delta > 0$, *there exists* $X = X(A, \delta)$ *such that*

$$\ell(y)/\ell(x) \le A\max((y/x)^\delta, (y/x)^{-\delta}) \qquad (x \ge X, y \ge X).$$

(ii) *If* $\ell$ *is bounded away from 0 and* $\infty$ *on every compact subset of* $[0, \infty[$, *then, for every* $\delta > 0$, *there exists* $A' = A'(\delta) > 1$ *such that*

$$\ell(y)/\ell(x) \le A'\max((y/x)^\delta, (y/x)^{-\delta}) \qquad (x \ge 0, y \ge 0).$$

**Acknowledgment.**   The authors are grateful to a referee for helpful comments.

A. GNEDIN
MATHEMATICAL INSTITUTE
UTRECHT UNIVERSITY
P.O. BOX 80010
3508 TA UTRECHT
THE NETHERLANDS
E-MAIL: gnedin@math.uu.nl

J. PITMAN
DEPARTMENT OF STATISTICS
UNIVERSITY OF CALIFORNIA
BERKELEY, CALIFORNIA 94720-3860
USA
E-MAIL: pitman@stat.berkeley.edu

M. YOR
LABORATORY OF PROBABILITY
UNIVERSITY OF PARIS VI
4 PLACE JUSSIEU, CASE 188
75252 PARIS, CEDEX 05
FRANCE
E-MAIL: deaproba@proba.jussieu.fr